\documentclass[12pt]{article}
\usepackage{latexsym}
\usepackage{amsmath}
\usepackage{amsfonts}
\usepackage{amsthm}
\usepackage{amssymb}

\def\0{\mathbf 0}
\def\1{\mathbf 1}

\numberwithin{equation}{section}
\newcounter{thm}[section]
\numberwithin{thm}{section}

\newtheorem{theorem}[thm]{Theorem}
\newtheorem{lemma}[thm]{Lemma}
\newtheorem{proposition}[thm]{Proposition}
\newtheorem{corollary}[thm]{Corollary}

\newtheorem{definition}[thm]{Definition}

 \title{Algebraic characterizations of measure algebras
 \thanks{Supported in part by the
 GAAV Grant IAA100190509}}
 \author{Thomas Jech \\
    Mathematical Institute, AS CR \\
    Zitna 25 \\
    CZ - 115 67 Praha 1 \\
    Czech Republic \\
    e-mail: jech@math.cas.cz
    }

\begin{document}

\maketitle

\abstract We present necessary and sufficient conditions for the
existence of a countably additive measure on a complete Boolean
algebra.
\endabstract

\section{Statement of results}\label{S1}

A \emph{Boolean algebra} is an algebra $B$ of subsets of a given
nonempty set $X$, with Boolean operations $a \cup b$, $a \cap b$,
$-a=X-a$, and the zero and unit elements $\mathbf 0=\emptyset$ and
$\mathbf 1=X$. A \emph{Boolean $\sigma$-algebra} is a Boolean
algebra $B$ such that every countable set $A\subset B$ has a
supremum $\sup A = \bigvee A$ (and an infimum $\inf A = \bigwedge
A$) in the partial ordering of $B$ by inclusion.

\begin{definition}\label{meas}
A \emph{measure} (more precisely, a strictly positive
$\sigma$-additive probability measure) on a Boolean
$\sigma$-algebra $B$ is a real valued function $m$ on $B$ such
that
\begin{itemize}
\item [(i)] $m(\mathbf 0) = 0,\ m(a) > 0 \text{ for } a \not = \mathbf 0
        \text{, and } m(\mathbf 1) = 1 $
\item [(ii)]$ m(a)\le m(b) \text{ if } a\subset b$
\item [(iii)] $m(a \cup b) =m(a) + m(b) \text{ if } a\cap b =\mathbf
0$
\item [(iv)]
$m(\bigvee_{n=1}^\infty a_n) = \sum_{n=1}^\infty m(a_n)
        \text{ whenever the } a_n \text{ are pairwise disjoint.}$
\end{itemize}

A \emph{measure algebra} is a Boolean $\sigma$-algebra that
carries a measure.
\end{definition}

Let $B$ be a Boolean algebra and let $B^+ = B - \{\mathbf 0\}$. A
set $A \subset B^+$ is an \emph{antichain} if $a\cap b =\mathbf 0$
whenever $a$ and $b$ are distinct elements of $A$. A
\emph{partition} $W$ (of $\mathbf 1$) is a maximal antichain, i.e.
an antichain with $\bigvee W =\mathbf 1$. $B$ satisfies the
\emph{countable chain condition} (ccc) if it has no uncountable
antichains. $B$ is \emph{weakly distributive} if for every
sequence $\{W_n\}_n$ of partitions there exists a partition $W$
with the property that each $a\in W$ meets only finitely many
elements of each $W_n$.

If $B$ is a measure algebra then $B$ satisfies ccc and is weakly
distributive. Below we present additional, purely algebraic,
conditions that characterize measure algebras.

If $\{a_n\}_n$ is a sequence in a Boolean $\sigma$-algebra $B$,
one defines $$ \limsup_n a_n = \bigwedge_{n=1}^\infty \bigvee_{k =
n}^\infty a_k \text{ and }\liminf_n a_n = \bigvee_{n=1}^\infty
\bigwedge_{k = n}^\infty a_k,$$ and if $\limsup_n a_n = \liminf_n
a_n =a$, then $a$ is the \emph{limit} of the sequence, denoted
$\lim_n a_n$.

\begin{theorem} \label{Thm1}
A Boolean $\sigma$-algebra $B$ is a measure algebra if and only if
it is weakly distributive and $B^+$ is the union of a countable
family $\{C_n\}_n$ such that for every $n$,
\begin{itemize}
\item [(i)] every antichain in $C_n$ has at most $K(n)$ elements
(for some integer $K(n)$), and
\item [(ii)] if $\{a_n\}_n$ is a sequence with $a_n \notin C_n$
for each $n$, then $\lim_n a_n =\0$.
\end{itemize}
\end{theorem}

\begin{theorem} \label{Thm5}
A Boolean $\sigma$-algebra $B$ is a measure algebra if and only if
$B^+$ is the union of a countable family $\{C_n\}_n$ such that for
every $n$,
\begin{itemize}
\item [(i)] every antichain in $C_n$ has at most $K(n)$ elements
(for some integer $K(n)$),
\item [(ii)] if $\{a_n\}_n$ is a sequence with $a_n \notin C_n$
for each $n$, then $\lim_n a_n =\0$, and
\item [(iii)] for every $k$, if $\{a_n\}_n$ is a sequence with
$\lim_n a_n =\0$, then for eventually all $n$, $a_n \notin C_k$.
\end{itemize}
\end{theorem}

\begin{theorem} \label{Thm2}
A Boolean $\sigma$-algebra $B$ is a measure algebra if and only if
it is weakly distributive and $B^+$ is the union of a countable
family $\{C_n\}_n$ such that for every $n$,
\begin{itemize}
\item [(i)] every antichain in $C_n$ has at most $K(n)$ elements,
and
\item [(ii)] for every $n$ and all $a$ and $b$, if $a\cup b\in
C_n$ then either $a\in C_{n+1}$ or $b\in C_{n+1}$
\end{itemize}
\end{theorem}

If a Boolean $\sigma$-algebra $B$ satisfies ccc then weak
distributivity is equivalent to this condition: if $\{W_n\}_n$ is
a sequence of partitions then each $W_n$ has a finite subset $E_n$
such that $\lim_n \bigcup E_n = \1$.

\begin{definition}
A Boolean $\sigma$-algebra $B$ is \emph{uniformly weakly
distributive} if there exists a sequence of functions $\{F_n\}_n$
such that for each partition $W$, $F_n(W)$ is a finite subset of
W, and if $\{W_n\}_n$ is a sequence of countable partitions then $\lim_n
\bigcup F_n(W_n) = \1$.
\end{definition}

\begin{definition} Let $B$ be a Boolean $\sigma$-algebra. $B$ is
\emph{concentrated} if for every sequence $A_n$ of finite
antichains with $|A_n| \ge 2^n$ there exist $a_n \in A_n$ such
that $\lim_n a_n = \0$. $B$ is \emph{uniformly concentrated} if
there exists a function $F$ such that for each finite antichain
$A$, $F(A)$ is an element of $A$, and if $A_n$ is a sequence of
finite antichains with $|A_n| \ge 2^n$ then $\lim_n F(A_n) = \0$.
\end{definition}

\begin{theorem} \label{Thm3}
A Boolean $\sigma$-algebra $B$ is a measure algebra if and only if
it is uniformly weakly distributive and concentrated.
\end{theorem}

\begin{theorem} \label{Thm4}
A Boolean $\sigma$-algebra $B$ is a measure algebra if and only if
it is weakly distributive and uniformly concentrated.
\end{theorem}

\section{Background and definitions}\label{S2}

We give a brief history of the problem, introduce relevant
definitions and state the known results. For an additional
reference and a more detailed history, see \cite{Handbook} (in
particular Fremlin's article \cite{Fr}) and \cite{BJ}.

The problem of an algebraic characterization of measure algebras
originated with John von Neumann. In 1937 (Problem 163 in
\cite{Sc}) he stated that measure algebras satisfy ccc and are
weakly distributive, and asked if these conditions are sufficient
for the existence of a measure.

In \cite{Mah}, Dorothy Maharam investigated Boolean
$\sigma$-algebras that carry a continuous submeasure and presented
necessary and sufficient conditions for the existence of such a
submeasure as well as of a measure.

\begin{definition}\label{submeas}
A (strictly positive) \emph{submeasure}  on a Boolean algebra $B$
is a real valued function $m$ on $B$ such that
\begin{itemize}
\item [(i)] $m(\mathbf 0) = 0,\ m(a) > 0 \text{ for } a \not = \mathbf 0
        \text{, and } m(\mathbf 1) = 1 $
\item [(ii)]$ m(a)\le m(b) \text{ if } a\subset b$
\item [(iii)] $m(a \cup b) \le m(a) + m(b)$
\end{itemize}
A \emph{Maharam submeasure} on a Boolean $\sigma$-algebra is a
submeasure that is continuous:
\begin{itemize}
\item [(iv)] if  $\{a_n\}_n$ is a decreasing sequence in $B$ with
$\bigwedge_{n=1}^\infty a_n = \0$ then $\lim_n m(a_n)=\0.$
\end{itemize}
A \emph{Maharam algebra} is a Boolean $\sigma$-algebra that
carries a Maharam submeasure.
\end{definition}

A measure is a Maharam submeasure, and every Maharam algebra is
ccc and weakly distributive. Maharam asked if every Maharam
algebra is a measure algebra. She also proved that a Suslin line,
if it exists, provides an example of a Boolean $\sigma$-algebra
that is ccc and weakly distributive but not a Maharam algebra.
(The existence of a Suslin line is consistent with the axioms of
set theory \cite{Ten}, \cite{Jech}, but not provable in ZFC,
\cite{ST}.)

In \cite{Ke}, John Kelley gave a combinatorial characterization of
Boolean algebras that carry a finitely additive measure. A
\emph{finitely additive measure} on a Boolean algebra is a
function $m$ that satisfies conditions (i), (ii) and (iii) of
Definition \ref{meas}. He also proved the following theorem (due
independently to Pinsker \cite{KVP}):

\begin{theorem}
A Boolean $\sigma$-algebra $B$ carries a measure if an only if it
is weakly distributive and carries a finitely additive measure.
\end{theorem}

\begin{proof}
Let $m$ be a finitely additive measure on $B$. For every partition
$W$ let $\lambda_W$ be the function $$ \lambda_W(b) = \sup
\{m(b\cap\bigcup E): E \text{ is a finite subset of }W\}, $$ and
let $$ \mu(b)=\inf\{\lambda_W(b): W \text{ is a partition}\}.$$
The function $\mu$ is $\sigma$-additive, and weak distributivity
implies that $\mu(b)>0$ for every $b\ne\0.$
\end{proof}

A major advance toward the solution of Maharam's problem was the
following result of Nigel Kalton and James Roberts.

\begin{definition}
A submeasure $m$ on a Boolean algebra $B$ is \emph{exhaustive} if
$\lim_n m(a_n) = 0$ for every infinite antichain $A = \{a_n :
n=1,2,...\}$. It is \emph{uniformly exhaustive} if for every
$\varepsilon > 0$ there exists some $n$ such that there is no
sequence of $n$ disjoint elements $a_1,\dots, a_n \in B$ with
$m(a_i) \geq \varepsilon$ for all $i = 1,\dots,n$.
\end{definition}

Note that a Maharam submeasure is exhaustive while a finitely
additive measure is uniformly exhaustive.

\begin{theorem}(Kalton-Roberts \cite{KR}.)
If a Boolean algebra $B$ carries a uniformly exhaustive submeasure
then $B$ carries a finitely additive measure.
\end{theorem}

\begin{corollary} If a Boolean $\sigma$-algebra carries a
uniformly exhaustive Maharam submeasure then it is a measure
algebra.
\end{corollary}

The use of the order sequential topology on $B$ (introduced by
Maharam and developed by Bohuslav Balcar) resulted in further
characterizations of Maharam algebras, cf. \cite{BGJ} and
\cite{BJP}, in particular:

\begin{theorem} (Balcar-Jech-Paz\'ak \cite{BJP}.) \label{B-J-P}
A Boolean $\sigma$-algebra $B$ carries a Maharam submeasure if an
only if it is ccc and weakly distributive and has the $G_\delta$
property, i.e. there exists a countable family $\{U_n\}_n$of
subsets of $B$ with $\bigcap_{n=1}^\infty U_n = \{\0\}$ such that
for every sequence $\{a_k\}_k$ with limit $\0$, eventually all
$a_k$ are in $U_n$.
\end{theorem}

Theorem \ref{B-J-P} combined with an earlier result of Todorcevic
\cite{Tod} shows that it is consistent that a Boolean
$\sigma$-algebra is a Maharam algebra if and only if it satisfies
ccc and is weakly distributive.

In \cite{Tal}, Michel Talagrand solved Maharam's problem by
constructing a submeasure on a countable Boolean algebra that is
exhaustive but not uniformly exhaustive. In view of \cite{KR} this
yields a (countably generated) Maharam algebra that is not a
measure algebra.

The Kalton-Roberts theorem and the Balcar-Jech-Paz\'ak theorem are
the tools we use in the proofs of Theorems \ref{Thm1}--\ref{Thm4}.

\section{Proof of Theorems \ref{Thm1}--\ref{Thm4}}\label{S3}

First we verify that measure algebras satisfy the conditions
stated in the theorems. Let $m$ be a measure on a Boolean
$\sigma$-algebra $B$. $B$ is weakly distributive, in fact
uniformly weakly distributive: For each $n$ and every partition
$W$, let $F_n(W)$ be a finite subset $E$ of $W$ such that
$m(\bigcup E)\ge 1-\frac{1}{2^n}$. If $\{W_n\}_n$ is a sequence of
partitions and if $a_n=-\bigcup F_n(W_n)$ then we have
$m(a_n)\le\frac{1}{2^n}$ and so $\limsup_n a_n =\0.$ Hence
$\lim_n\bigcup F_n(W_n)=\1.$

For each $n$ let $C_n$ be the family of all $a\in B$ such that
$m(a)\ge\frac{1}{2^n}.$ We have $\bigcup_{n=1}^\infty C_n = B^+$,
for every $n$ every antichain in $C_n$ has at most $2^n$ elements,
and if $a\cup b\in C_n$ then either $a\in C_{n+1}$ or $b\in
C_{n+1}.$ If $a_n\notin C_n$ for every $n$ then $\lim_n a_n=\0.$

For every finite antichain $A$ let $F(A)=a\in A$ be such that
$m(a)\le m(x)$ for all $x\in A.$ We have
$m(F(A))\le\frac{1}{|A|}$, and so if $\{A_n\}_n$ is a sequence of
finite antichains with $|A_n|\ge 2^n$ then for each $n$,
$m(F(A_n))\le\frac{1}{|2^n|}$ and it follows that $\lim_n
F(A_n)=\0.$

We shall prove that the conditions in Theorems
\ref{Thm1}--\ref{Thm4} imply the existence of a measure.

\begin{lemma} \label{unif}
Let $B$ be a weakly distributive Boolean $\sigma$-algebra that
satisfies the conditions of Theorem \ref{Thm1}. Then $B$ is
uniformly weakly distributive.
\end{lemma}

\begin{proof}
Let $\{C_n\}_n$ be a countable family that has properties (i) and
(ii). Without loss of generality we may assume that each $C_n$ is
upward closed, i.e. if $a\subset b$ and $a\in C_n$ then $b\in
C_n.$ To begin with, condition (i) implies ccc and so every
antichain is at most countable. Let $W$ be a partition and $n$ a
number. We shall define $F_n(W)$ so that the functions $F_n$
witness uniform weak distributivity.

We claim that there exists a finite set $E\subset W$ (possibly
empty) such that there exists no nonempty finite set $F\subset
W-E$ with $\bigcup F\in C_n$. If not then we can find an infinite
sequence $\{E_k\}_k$ of disjoint finite subsets of $W$ producing
an infinite antichain $\{\bigcup E_k: k=1,2,...\}$ in $C_n.$ We
let $F_n(W)$ be such an $E.$

Now let $\{W_n\}_n$ be a sequence of partitions. Since $B$ is
weakly distributive there exist finite sets $E_n\subset W_n$ such
that $\lim_n \bigcup E_n = \1.$ For each $n$ let $a_n=\bigcup E_n
- \bigcup F_n(W).$ By the definition of $F_n(W)$ we have
$a_n\notin C_n$ and hence $\lim_n a_n =\0.$ It follows that
$\lim_n \bigcup F_n(W)=\1.$
\end{proof}

\begin{lemma}(\cite{BJ}, p. 259.) \label{Gdelta}
If $B$ is a uniformly weakly
distributive ccc Boolean $\sigma$-algebra then $B$ has the
$G_\delta$ property.
\end{lemma}

\begin{proof}
Let $F_n$ be functions that witness the uniform weak
distributivity. For each $n$ we let $$U_n=\{a\in B: a \text{ is
disjoint from } \bigcup F_n(W) \text{ for some partition } W\}.$$
First we claim that $\bigcap_{n=1}^\infty U_n = \{\0\}:$ If $a\in
U_n$ for each $n,$ and if $W_n$ are partitions such that $a\cap
F_n(W_n)=\0$ then because $\lim_n F_n(W_n)=\1$, $a$ must be $\0.$
Now let $\{a_k\}_k$ be a sequence with limit $\0,$ and let $n$ be
an integer. There is a decreasing sequence $\{b_k\}_k$ such that
$b_1=\1,$ $a_k\subset b_k$ for each $k$ and
$\bigwedge_{k=1}^\infty b_k=\0.$ Let $W$ be the partition
$\{b_{k+1}-b_k: k=1,2,...\}$ and let $E=F_n(W).$ There is some $K$
such that $b_K\cap\bigcup E =\0$ and hence $a_k\in U_n$ for all
$k\ge K.$
\end{proof}

\begin{lemma} \label{unifexh}
Let $B$ be a Maharam algebra with a Maharam submeasure $m$ and
assume that $B$ satisfies the conditions of Theorem \ref{Thm1}.
Then $m$ is uniformly exhaustive.
\end{lemma}

\begin{proof}
Let $\{C_n\}_n$ be a countable family with properties (i) and
(ii). In order to verify that $m$ is uniformly exhaustive it
suffices to show that for every $\varepsilon>0$ there is some $n$
such that $\{a\in B: m(a)\ge\varepsilon\} \subset C_n.$ If not,
let $\varepsilon$ be a counterexample. For each $n$ we pick
$a_n\notin C_n$ with $m(a_n)\ge\varepsilon.$ By (ii), $\lim_n a_n
=\0.$ Since $m$ is continuous, we have $\lim_n m(a_n)=0$, a
contradiction.
\end{proof}

Now Theorem \ref{Thm1} follows: If $B$ satisfies the conditions,
then by Lemmas \ref{unif}, \ref{Gdelta} and the
Balcar-Jech-Paz\'ak Theorem $B$ carries a Maharam submeasure, and
by Lemma \ref{unifexh} and the Kalton-Roberts Theorem, $B$ carries
a measure.

\bigskip

Having proved Theorem \ref{Thm1}, for Theorem \ref{Thm5} it
suffices to show that under the conditions of the theorem, $B$ is
weakly distributive. As $B$ satisfies ccc it is enough to show
that if for every $k$, $\lim_n a_n^k =\0$, then there is a
function $n(k)$ such that $\lim_k a_{n(k)}^k = \0$ (see \cite{BJ},
p.253). This ``diagonal property'' is verified using (iii) and
(ii).

\bigskip

Turning our attention to Theorem \ref{Thm2}, we will show that the
conditions of Theorem \ref{Thm2} imply the conditions of Theorem
\ref{Thm1}. Let $B$ be a weakly distributive Boolean
$\sigma$-algebra and let $\{C_n\}_n$ be a countable family that
has properties (i) and (ii) of Theorem \ref{Thm2}. Notice that if
we replace each $C_n$ by the set $\{x\in B^+: (\exists y \subset
x)\, y\in C_1\cup\dots\cup C_n\},$ then the family still has
properties (i) and (ii). Thus we assume that $C_1\subset
C_2\subset\dots$ and that each $C_n$ is upward closed. The
following lemma shows that $\{C_n\}_n$ satisfies condition (ii) of
Theorem \ref{Thm1}.

\begin{lemma}
If $a_n\notin C_n$ for each $n$, then $\limsup_n a_n=\0.$
\end{lemma}

\begin{proof}
Let $a=\limsup_n a_n$ and assume that $a\ne \0.$ For each $n$ and
each $k$, let $b_{nk}=a_{n+1}\vee\dots\vee a_{n+k}.$ From (ii) it
follows that $b_{nk}\notin C_n$, for all $k.$

We have $a=\lim_n \lim_k b_{nk}$, and by weak distributivity there
exists for each $n$ some $k(n)$ such that $a=\lim_n b_{n,k(n)}.$
Since $a\ne\0$, there exist some $b\ne\0$, $b\subset a$ and some
$N$ such that $b\subset b_{n,k(n)}$ for all $n\ge N.$ Let $n\ge N$
be such that $b\in C_n.$ Since $C_n$ is upward closed, we have
$b_{n,k(n)}\in C_n$, a contradiction.
\end{proof}

\bigskip

For Theorem \ref{Thm3}, let $B$ be a Boolean $\sigma$-algebra that
satisfies the conditions of Theorem \ref{Thm3}. By Lemma
\ref{Gdelta} $B$ has the $G_\delta$ property, and Theorem
\ref{B-J-P} shows that $B$ is a Maharam algebra as long as it is
ccc. We use the following lemma:

\begin{lemma} (\cite{BJ2})\label{ccc}
If $B$ is a uniformly weakly
distributive Boolean $\sigma$-algebra then $B$ satisfies ccc.
\end{lemma}

\begin{proof} Let $\bar B$ be the regular completion of $B$. Since
$B$ is dense in $\bar B$, every partition in $\bar B$ has a
refinement in $B$ and limits of sequences in $B$ are the same in
$\bar B$ as in $B$. Hence $\bar B$ is uniformly weakly
distributive. If $\bar B$ has a partition of size $\omega_1$ then
$P(\omega_1)$ is a complete subalgebra of $\bar B$ and therefore
it is uniformly weakly distributive. By \cite{BJ2} $P(\omega_1)$
is not uniformly weakly distributive and so $\bar B$, and hence
$B$, satisfies ccc.
\end{proof}

Hence $B$ carries a Maharam submeasure by the Balcar-Jech-Paz\'ak
Theorem.

\begin{lemma} \label{unifexh2}
Let $B$ be a Maharam algebra with a Maharam submeasure $m$ and
assume that $B$ is concentrated. Then $m$ is uniformly exhaustive.
\end{lemma}

\begin{proof}
If $m$ is not uniformly exhaustive then there exists an
$\varepsilon>0$ such that for every $n$ there is an antichain
$A_n$ of size $2^n$ with $m(a)\ge\varepsilon$ for all $a\in A_n$.
This contradicts the condition that there exists a sequence
$\{a_n\}_n$ such that $a_n\in A_n$ and $\lim_n a_n=\0.$
\end{proof}

This completes the proof of Theorem \ref{Thm3}: if $B$ satisfies
its conditions then $B$ carries a Maharam submeasure by the
Balcar-Jech-Paz\'ak Theorem, and by the Kalton-Roberts Theorem it
carries a measure.

\bigskip

For Theorem \ref{Thm4}, we will show that if $B$ satisfies the
conditions of Theorem \ref{Thm4} then it satisfies the conditions
of Theorem \ref{Thm1}. Let $B$ be a weakly distributive Boolean
$\sigma$-algebra and let $F$ be a function acting on finite
antichains  witnessing that $B$ is uniformly concentrated.

For each $n$ we let $$ C_n=\{a\in B^+: a\ne F(A) \text{ for every
antichain } A \text { of size } \ge 2^n\}.$$ If $a$ is such that
$a\notin C_n$ for each $n$ then there exist antichains $A_n$ such
that $|A_n|\ge 2^n$ and $a=F(A_n).$ Since $\lim_n F(A_n)=\0$ we
have $a=\0$, and so $\bigcup_{n=1}^\infty C_n=B^+.$ If $\{a_n\}_n$
is a sequence such that $a_n\notin C_n$ for each $n$ then there
exist antichains $A_n$ such that $|A_n|\ge 2^n$ and $a_n=F(A_n).$
Hence $\lim_n a_n=\0$. Finally, every antichain in $C_n$ has fewer
than $2^n$ elements: If $A$ is an antichain of size $\ge 2^n$,
then $F(A)\notin C_n$ and so $A$ is not a subset of $C_n$. Hence
$B$ satisfies the assumptions of Theorem \ref{Thm1}.

\section{Odds and ends}

A Boolean algebra $B$ satisfies the \emph{$\sigma$-bounded cc
(chain condition)} if $B^+$ is the union of a countable family
$\{C_n\}_n$ such that for every $n$, every antichain in $C_n$ has
at most $K(n)$ elements (for some integer $K(n)$). $B$ satisfies
the \emph{$\sigma$-finite cc} if $B^+$ is the union of a countable
family $\{C_n\}_n$ such that for every $n$, every antichain in
$C_n$ if finite. These conditions were explicitly stated in
\cite{HT}.

The conditions in Theorems \ref{Thm1} and \ref{Thm2} state that
$B$ is $\sigma$-bounded cc but require that the $C_n$ have an
additional property. This is necessary: Talagrand's result
\cite{Tal} yields a Maharam algebra that is $\sigma$-bounded cc
but is not a measure algebra. In contrast, Stevo Todorcevic proved
in \cite{Tod2} that a Boolean $\sigma$-algebra $B$ is a Maharam
algebra if and only if it is weakly distributive and
$\sigma$-finite cc. Notice that if $B$ carries a Maharam
submeasure $m$ then $B^+=\bigcup_{n=1}^\infty C_n$ such that
$\{C_n\}_n$ witnesses the $\sigma$-finite cc and also has the
additional properties from Theorems \ref{Thm1} and \ref{Thm2}.
(Let $C_n=\{a\in B: m(a)\ge \frac{1}{2^n}\}.$)

As David Fremlin pointed out, if we drop weak distributivity in
Theorem \ref{Thm2}, then we get a characterization of Boolean
algebras that carry a finitely additive measure:

\begin{theorem}
A Boolean algebra $B$ carries a finitely additive measure if and
only if  $B^+$ is the union of a countable family $\{C_n\}_n$ such
that for every $n$,
\begin{itemize}
\item [(i)] every antichain in $C_n$ has at most $K(n)$ elements,
and
\item [(ii)] for every $n$ and all $a$ and $b$, if $a\cup b\in
C_n$ then either $a\in C_{n+1}$ or $b\in C_{n+1}$
\end{itemize}
\end{theorem}

\begin{proof}
The condition is clearly necessary. For the sufficiency we follow
the construction in \cite{BGJ}, p. 75: First, modify the $C_n$ so
that each $C_n$ is upward closed and $C_1\subset C_2\subset C_3
\subset \dots$, and let $U_n=B-C_n$ for each $n$. For each number
$r=\frac{1}{2^{n_1}}+\dots+\frac{1}{2^{n_k}}$ iwth
$n_1<\dots<n_k,$ let $V_r=\{u_1\cup\dots\cup u_k: u_i\in
U_{n_i},i=1,\dots k\},$ and define $m(a)=\inf\{r: a\in V_r\}.$ The
condition (ii) implies that for each $a$, $\frac{1}{2^n}\le
m(a)\le \frac{1}{2^{n-1}}$, where $n$ is the least $n$ such that
$a\in C_n.$ It follows that $m(a)>0$ whenever $a\ne\0,$ and $m$ is
a submeasure on $B$. By (i), $m$ is uniformly exhaustive, and so
by the Kalton-Roberts Theorem $B$ carries a finitely additive
measure.
\end{proof}

Theorems \ref{Thm3} and \ref{Thm4} state that measure algebras are
characterized by being uniformly weakly distributive and
concentrated, resp. weakly distributive and uniformly
concentrated. From Theorem \ref{Thm3} and \cite{BJP} it follows
that it is consistent that $B$ is a measure algebra if and only if
it is ccc, weakly distributive and concentrated. On the other
hand, if a Suslin tree exists then the corresponding Suslin
algebra $B$ is ccc, weakly distributive and concentrated, and does
not carry even a Maharam submeasure:

\begin{proposition} Let $T$ be a Suslin tree and $B$ the
corresponding complete Boolean algebra. Then $B$ is concentrated.
\end{proposition}

\begin{proof}
Let $\{A_n\}_n$ be finite antichains in $B$, $|A_n|\ge 2^n.$ We
may assume that $\bigcup A_n = \1$ for each $n$. A routine
argument using that $T$ is a Suslin tree shows that there exists a
countable family of functions $\{f_k\}_k \subset
\prod_{n=1}^\infty A_n$ such that $$ \bigvee_{k=1}^\infty
\bigwedge_{n=1}^\infty f_k(n) = \1.$$ Now let $F
\in\prod_{n=1}^\infty A_n$ be such that for each $n$, $F(n)\ne
f_i(n)$ for all $i=1,\dots, n$. We show that $\lim_n F(n)=\0.$

Let $a_n=F(n)$ and let $a=\limsup_n a_n.$ For each $k$ let
$b_k=\bigwedge_{n=1}^\infty f_k(n).$ Since $a_n\cap f_k(n)=\0$ for
all $n\ge k$, we have $\bigvee_{n=k}^\infty a_n \cap b_k = \0$,
and it follows that $a\cap b_k=\0.$ Hence $a=\0.$
\end{proof}

In the proof of Theorem \ref{Thm4} we showed that if $B$ is
uniformly concentrated then $B$ is $\sigma$-bounded cc. It turns
out that a weak version of uniformly concentrated is equivalent to
the $\sigma$-bounded cc, and uniformly concentrated is equivalent
to conditions (i) and (ii) of Theorem \ref{Thm1}:

\begin{proposition}
(a) A Boolean algebra $B$ is $\sigma$-bounded cc if and only if
there exists a function $F$ such that for each finite antichain
$A$, $F(A)\in A$, and if $\{A_n\}_n$ is a sequence of finite
antichains of increasing size then $\liminf_n F(A_n)=\0.$

(b) A Boolean algebra $B$ is uniformly concentrated if and only if
$B^+$ is the union of a countable family $\{C_n\}_n$ with
$C_1\subset C_2\subset\dots$ such that for every $n$,
\begin{itemize}
\item [(i)] every antichain in $C_n$ has fewer than $2^n$
elements, and
\item [(ii)] if $\{a_n\}_n$ is a sequence with $a_n \notin C_n$
for each $n$, then $\lim_n a_n =\0$.
\end{itemize}
\end{proposition}

\begin{proof}
(a) First assume that $B$ is $\sigma$-bounded cc, and let
$\{C_n\}_n$ be a witness. We may assume that each $C_n$ is upward
closed and that $C_1\subset C_2\subset...$ If $A$ is a finite
antichain, let $n$ be the least $n$ such that $A\subset C_n$ and
let $F(A)=a\in A$ be such that $a\notin C_{n-1}$. Hence for all
$N$, if $F(A)\in C_N$ then $A\subset C_N.$

Now let $\{A_n\}_n$ be a sequence of finite antichains increasing
in size, and let $a_n=F(A_n).$ We claim that $\liminf_n a_n=\0.$
If not then there exist some $a\ne \0$ and some $k$ such that
$a\subset a_n$ for all $n\ge k.$ Let $N$ be such that $a\in C_N$;
then $a_n\in C_N$ for all $n\ge k.$ It follows that $A_n\subset
C_N$ for all $n\ge k$, and so $C_N$ has antichains of arbitrary
size, a contradiction.

Conversely, let $F$ be a function that satisfies the condition. If
we let $C_n=\{a\in B^+: a\ne F(A)\text{ for every antichain of
size } \ge n+1\}$ then the same argument we used in the proof of
Theorem \ref{Thm4} shows that $\bigcup_{n=1}^\infty C_n = B^+$ and
that every antichain in $C_n$ has at most $n$ elements.

(b) For one direction, see the proof of Theorem \ref{Thm4}. For
the other direction, given the $C_n$, we let $F(A)=a\in A$ be such
that $a\notin C_{n-1}$ where $n$ is the least $n$ with $A\subset
C_n.$ Now if $|A_n|\ge 2^n$ then $A_n\not\subset C_n$ and so
$F(A_n)\notin C_n.$ Hence  $\lim_n F(A_n)=\0.$
\end{proof}

Weak distributivity has a formulation in terms of forcing: a
complete ccc Boolean algebra $B$ is weakly distributive if and
only if for every $B$-name $\dot{f}$ for a function from $\omega$
to $\omega$ there exists a function $g:\omega \rightarrow \omega$
such that $$\Vdash \exists N\forall n\ge N \dot{f}(n) < g(n).$$
(The last formula is equivalent to $\lim_n
||\dot{f}(n)<g(n)||=\1.$)

Similarly, $B$ is concentrated if and only if for every $B$-name
$\dot{f}$ for a function from $\omega$ to $\omega$ there exists a
function $g:\omega \rightarrow \omega$ such that $g(n)<2^n$ for
each $n$ and $$\Vdash \exists N\forall n\ge N \dot{f}(n)\ne
g(n).$$

The following result shows that the existence of a finitely
additive measure does not imply that $B$ is concentrated. The
Cohen algebra carries a finitely additive measure but is not
concentrated:

\begin{proposition}
The Cohen algebra is not concentrated.
\end{proposition}

\begin{proof}
We use this representation of the Cohen algebra: Let $P$ be the
forcing where the forcing conditions are finite sequences $p$ of
integers such that $p(n)<2^n$ for each $n\in dom(p).$ We let
$\dot{f}$ be the following name for a function from $\omega$ to
$\omega$: for each $n$ and each $k<2^n$ let $$ ||\dot{f}(n)=k|| =
\bigvee \{p:p(n)=k\}.$$ Now if $g:\omega \rightarrow \omega$ is
such that $g(n)<2^n$ for all $n$ then for every condition $p$ and
every $N$ there exist a stronger condition $q$ and some $n>N$ such
that $q \Vdash \dot{f}(n)=g(n).$ This shows that $\dot{f}$ is a
counterexample.
\end{proof}

Maharam algebras have a characterization in terms of infinite
games. Using \cite{BGJ}, David Fremlin proved in \cite{Fr2} that a
strategic version of weak distributivity implies the existence of
a Maharam submeasure for Boolean $\sigma$-algebras that satisfy
ccc (see \cite{BJ}, p. 261, for details). In \cite{BJ2} it is
shown that the ``strategic diagonal property'' implies ccc.
Combining this with the proof of Theorem \ref{Thm3}, we obtain the
following characterization of measure algebras:

Let $B$ be a complete Boolean algebra and consider the infinite
game $\mathcal G$ in which the $n$th move of Player I is a
$B$-name $\dot{f}(n)$ for an integer and the $n$th move of Player
II is an integer $g(n)$. Thus I produces a $B$-name $\dot{f}$ for
a function from $\omega$ to $\omega$ and II produces a function
$g:\omega \rightarrow \omega$. Player II wins if $$\Vdash \exists
N\forall n\ge N\, \dot{f}(n)< g(n) \text{ and }
(\dot{f}(n)\not\equiv g(n) \mod{2^n}).$$

\begin{theorem} A complete Boolean algebra $B$ is a measure
algebra if and only if Player II has a winning strategy in the
game $\mathcal G$.
\end{theorem}

\bibliographystyle{plain}

\end{document}